\documentclass[10pt,reqno]{amsart}
\usepackage{amsfonts,amsmath}
\usepackage{graphicx}
\usepackage{subfigure}
\usepackage{lscape}

\begin{document}

\begin{center}
{\LARGE Putting the Prisoner's Dilemma in Context \vspace{0.7cm}}
\end{center}

{\large L. A. Khodarinova}$^{{\dag }}${\large \ and\ J. N. Webb
{\large \ \vspace{0.3cm}}

$^{{\dag }}$\textit{Magnetic Resonance Centre, School of Physics and
Astronomy, University of Nottingham, Nottingham, England NG7 2RD, e-mail:
LarisaKhodarinova@hotmail.com}
\vspace{1.0cm}

\textbf{Abstract.} The standard iterated prisoner's dilemma is an unrealistic model of social
behaviour because it forces individuals to participate in the interaction.
We analyse a model in which players have the option of ending their
association. If the payoff for living alone is neither too high nor too low
then the potential for cooperative behaviour is enhanced. For some parameter
values it is also possible for a polymorphic population of defectors and
conditional cooperators to be stable.

\section*{Introduction}

The iterated, or repeated, prisoner's dilemma is the most popular model of
social interactions \cite{AH}. Since its inception the basic model has been
modified in many ways (see Dugatkin \cite{D} for a review). However, in all
these versions it is assumed that the players must engage in the interaction
and have no opportunity to end it. This unrealistic feature is just one
facet of the more general assumption that one particular social interaction
may be considered in isolation from all others that an individual may face.

In this paper we use the framework of stochastic games \cite{FV} to consider
a version of the iterated prisoner's dilemma in which the players may choose
to discontinue their association. We assume that once the partnership has
been dissolved by one or more of the players, then each receives the same,
fixed, per-period payoff. This is probably the simplest way that an
interaction can be considered as being dependent on other situations in
which individuals find themselves during a complex and, at least partly,
social life. We will use the standard replicator dynamics \cite{TJ}\ to
investigate the effect that the existence of this outside option has on the
evolution of cooperative behaviour in a population of players.

\section*{The Model}

A general stochastic game has three major components: the set of states, the
games played in each of these states and the (possibly behaviour-dependent)
probabilities for transition between the states. In our model the states
represent the different contexts in which players may interact, so we will
refer to them as context games.

We consider an interaction described by the following multi-state,
stochastic game. There are three possible context-games (states) $G_{0},$ $%
G_{1}$ and $G_{2}.$ The interaction starts with context-game $G_{0}.$ In
this game the players make the decision about whether or not they wish to
initiate or continue an association. The first player and the second player
choose between two possible actions: $A$=\textquotedblleft
associate\textquotedblright\ or $B$=\textquotedblleft break
up\textquotedblright . There are no payoffs directly associated with this
decision. Context-game $G_{1}$ represents some specific activity in which
the individuals can participate together. It is modelled by the prisoner's
dilemma and the players choose between the possible actions: $C$%
=\textquotedblleft cooperation\textquotedblright\ or $D$=\textquotedblleft
defection\textquotedblright . Context-game $G_{2}$ can be considered as a
background state representing the situation when there is no interaction or
association between the players. There is only one possible action: $L$%
=\textquotedblleft be alone\textquotedblright. \medskip

\begin{center}
\begin{tabular}{c}
$%
\begin{array}{ll}
G_{0}\!: & 
\begin{tabular}{|l|c|c|}
\hline
{\footnotesize P}$_{1}\setminus ${\footnotesize P}$_{2}$ & {\footnotesize A}
& {\footnotesize B} \\ \hline
\multicolumn{1}{|c|}{\footnotesize A} & \multicolumn{1}{|l|}{$\!\!\!%
\rule[-0.14in]{0in}{0.35in}%
\begin{array}{l}
\!\!^{\left( 0,0\right) }\!\!\left/ \!_{\rule[-0.06in]{0in}{0.14in}\left[
0,1,0\right] }\!\!\right. \!%
\end{array}%
\!\!\!$} & \multicolumn{1}{|l|}{$\!\!\!\rule[-0.14in]{0in}{0.35in}%
\begin{array}{l}
\!\!^{\left( 0,0\right) }\!\!\left/ \!_{\rule[-0.06in]{0in}{0.14in}\left[
0,0,1\right] }\!\!\right. \!%
\end{array}%
\!\!\!$} \\ \hline
\multicolumn{1}{|c|}{\footnotesize B} & \multicolumn{1}{|l|}{$\!\!\!%
\rule[-0.14in]{0in}{0.35in}%
\begin{array}{l}
\!\!^{\left( 0,0\right) }\!\!\left/ \!_{\rule[-0.06in]{0in}{0.14in}\left[
0,0,1\right] }\!\!\right. \!%
\end{array}%
\!\!\!$} & \multicolumn{1}{|l|}{$\!\!\!\rule[-0.14in]{0in}{0.35in}%
\begin{array}{l}
\!\!^{\left( 0,0\right) }\!\!\left/ \!_{\rule[-0.06in]{0in}{0.14in}\left[
0,0,1\right] }\!\!\right. \!%
\end{array}%
\!\!\!$} \\ \hline
\end{tabular}%
\end{array}%
\begin{array}{ll}
G_{1}\!: & 
\begin{tabular}{|l|c|c|}
\hline
{\footnotesize P}$_{1}\setminus ${\footnotesize P}$_{2}$ & {\footnotesize C}
& {\footnotesize D} \\ \hline
\multicolumn{1}{|c|}{\footnotesize C} & \multicolumn{1}{|l|}{$\!\!\!%
\rule[-0.14in]{0in}{0.35in}%
\begin{array}{l}
\!\!^{\left( r,r\right) }\!\!\left/ \!_{\rule[-0.06in]{0in}{0.14in}\left[
1,0,0\right] }\!\!\right. \!%
\end{array}%
\!\!\!$} & \multicolumn{1}{|l|}{$\!\!\!\rule[-0.14in]{0in}{0.35in}%
\begin{array}{l}
\!\!^{\left( s,t\right) }\!\!\left/ \!_{\rule[-0.06in]{0in}{0.14in}\left[
1,0,0\right] }\!\!\right. \!%
\end{array}%
\!\!\!$} \\ \hline
\multicolumn{1}{|c|}{\footnotesize D} & \multicolumn{1}{|l|}{$\!\!\!%
\rule[-0.14in]{0in}{0.35in}%
\begin{array}{l}
\!\!^{\left( t,s\right) }\!\!\left/ \!_{\rule[-0.06in]{0in}{0.14in}\left[
1,0,0\right] }\!\!\right. \!%
\end{array}%
\!\!\!$} & \multicolumn{1}{|l|}{$\!\!\!\rule[-0.14in]{0in}{0.35in}%
\begin{array}{l}
\!\!^{\left( p,p\right) }\!\!\left/ \!_{\rule[-0.06in]{0in}{0.14in}\left[
1,0,0\right] }\!\!\right. \!%
\end{array}%
\!\!\!$} \\ \hline
\end{tabular}%
\end{array}%
$ \\ 
$%
\begin{array}{ll}
G_{2}\!: & \rule[0.2in]{0in}{0.35in}%
\begin{tabular}{|l|c|}
\hline
{\footnotesize P}$_{1}\setminus ${\footnotesize P}$_{2}$ & {\footnotesize L}
\\ \hline
\multicolumn{1}{|c|}{\footnotesize L} & \multicolumn{1}{|l|}{$\!\!\!%
\rule[-0.14in]{0in}{0.35in}%
\begin{array}{l}
\!\!^{\left( z,z\right) }\!\!\left/ \!_{\rule[-0.06in]{0in}{0.14in}\left[
0,0,1\right] }\!\!\right. \!%
\end{array}%
\!\!\!$} \\ \hline
\end{tabular}%
\ \medskip 
\end{array}%
$ \\ 
$%
\begin{array}{c}
\text{\rule{0in}{0.35in}Table1.The multi-state game.} \\ 
\text{In each state the per-period payoffs for the relevant action choices }
\\ 
\text{are given above and to the left of the diagonal line;} \\ 
\text{the behaviour-dependent probabilities of transition to the other
states } \\ 
\text{are given below the line.\rule[-0.3in]{0in}{0.18in}}%
\end{array}%
$%
\end{tabular}
\end{center}

The actions chosen define both immediate payoffs to the individuals and\
future transition probabilities. The immediate payoffs collected by the
players are given in table~1. The first entry in each payoff pair contains
the payoff to the player $P_{1},$ who selects the row action, the second is
for the player $P_{2},$ who selects the column action. In this paper we are
considering an extension of the standard iterated prisoner's dilemma, for
which the following inequalities hold in $G_{1}$. 
\begin{equation*}
t>r>p>s\geq 0.
\end{equation*}%
Transition probabilities, which are determined by the choice of actions are
presented in table~1 as a set of three numbers. This set of numbers appears
in square brackets in each cell of the matrices. Here the first, second or
third number is, respectively, the probability that context-game $G_{0},$ $%
G_{1}$ or $G_{2}$ is played at the next round. The probabilities are defined
by the following rules. If context-game $G_{0}$ is played and action $A=$%
\textquotedblleft associate\textquotedblright\ is chosen by both players, at
the next round context-game $G_{1}$ is played; if action $B$%
=\textquotedblleft break up\textquotedblright\ is chosen by at least one
player, context-game $G_{2}$ is played at the next round with probability~$%
1. $ Whatever actions are chosen when context-game $G_{1}\ $is played,
context-game $G_{0}$ is played at the next round with probability~$1$. If
context-game $G_{2}$ is played, at the next round context-game $G_{2}$ is
played again with probability$~1.$

We assume that after playing context game $G_{0}$, players survive to play
game $G_{1}$or $G_{2}$ (as appropriate) with probability 1. After playing
context games $G_{1}$ or $G_{2}$ players survive to the next round with
probability $\beta $ $(0\leq \beta <1)$. In principle, these survival
probabilities could be different but, for simplicity, we will assume they
are equal.

As with the iterated prisoner's dilemma there is an infinite number of pure
strategies that could be considered. We will initially restrict our
attention to the following three strategies.

\begin{itemize}
\item Conditional cooperation (which we denote $\sigma _{C}).$ A player
following this strategy will initially ``Associate'' in $G_{0}$ then
``Cooperate'' in $G_{1}$; if this behaviour is reciprocated then the player
will continue to associate and cooperate; otherwise\ it will choose ``Break
up'' in $G_{0}$.

\item Defection (which we denote $\sigma _{D}$). A player following this
rather pathological strategy will ``Associate'' in $G_{0}$ and then
``Defect'' in $G_{1}$.

\item An unsociable strategy (which we denote $\sigma _{B}).$ A player
following this strategy will ``Break up'' in $G_{0}$. Strictly speaking this
is a set of strategies since any behaviour is allowed in $G_{1}$. However,
since we do not consider the possibility that players make errors, the
behaviour in $G_{1}$ does not affect payoffs. Consequently we ignore this
technicality.
\end{itemize}

\noindent The consequences of introducing a fourth strategy of unconditional
cooperation will be considered later.

\section*{Evolutionary Dynamics}

We set up the evolutionary dynamics by considering an infinitely large
population of individuals who adopt one of the three pure strategies. The
payoffs in the repeated game, $\pi (\sigma ,\sigma ^{\prime })$\ for
adopting strategy $\sigma $ against an opponent who adopts strategy $\sigma
^{\prime }$ are 
\begin{equation}
A=\left[ 
\begin{array}{lll}
\pi (\sigma _{C},\sigma _{C}) & \pi (\sigma _{C},\sigma _{D}) & \pi (\sigma
_{C},\sigma _{B}) \\ 
\pi (\sigma _{D},\sigma _{C}) & \pi (\sigma _{D},\sigma _{D}) & \pi (\sigma
_{D},\sigma _{B}) \\ 
\pi (\sigma _{B},\sigma _{C}) & \pi (\sigma _{B},\sigma _{D}) & \pi (\sigma
_{B},\sigma _{B})%
\end{array}%
\right] =\frac{1}{1-\beta }\left[ 
\begin{array}{ccc}
r & s\left( 1-\beta \right) +\beta z & z \\ 
t\left( 1-\beta \right) +\beta z & p & z \\ 
z & z & z%
\end{array}%
\right] .  \label{Payoff_Matrix}
\end{equation}

Let $x_{1}$ and $x_{2}$ be the proportions of individuals who adopt $\sigma
_{C}$ and $\sigma _{D}$ respectively. The proportion of individuals using $%
\sigma _{B}$ is then $1-x_{1}-x_{2}$. The standard replicator dynamics \cite%
{TJ} is then two equations describing the evolution of a point $x=\left(
x_{1},x_{2}\right) $ in the domain 
\begin{equation}
\Delta =\left\{ \left( x_{1},x_{2}\right) :(x_{1}\geq 0)\cap (x_{2}\geq
0)\cap (x_{1}+x_{2}\leq 1)\right\} .  \label{Simplex_2}
\end{equation}

Denote 
\begin{eqnarray*}
a &=&\left( z-r\right) -\gamma \left( z-t\right) ;\quad b=\gamma \left(
z-s\right) -\left( z-p\right) ; \\
c &=&\left( z-r\right) ;\quad f=\gamma \left( z-s\right) \quad \text{where}%
\quad \gamma =1-\beta .
\end{eqnarray*}%
Then the Replicator Dynamics can be written as the following system of
equations. 
\begin{equation*}
\begin{array}{l}
\dot{x}_{1}=x_{1}\left( cx_{1}^{2}+\left( f+c-a\right) x_{1}x_{2}+\left(
f-b\right) x_{2}^{2}-cx_{1}-fx_{2}\right)  \\ 
\dot{x}_{2}=x_{2}\left( cx_{1}^{2}+\left( f+c-a\right) x_{1}x_{2}+\left(
f-b\right) x_{2}^{2}+\left( a-c\right) x_{1}+\left( b-f\right) x_{2}\right) 
\end{array}%
\end{equation*}%
Although this system is integrable for arbitrary choices of parameter values 
\cite{B}, the general solution given in appendix A is not easy to work with.
We will now introduce a commonly used set of values for the prisoner's
dilemma context game $G_{1}$ to reduce the number of parameters, and we will
use the standard linearization approach to study how the solution depends on
the value of the outside option, $z$, and the survival probability, $\beta $%
. Accordingly we set 
\begin{equation*}
r=3,\qquad s=0,\qquad t=5,\qquad p=1.
\end{equation*}%
The payoff matrix then becomes 
\begin{equation*}
A=\frac{1}{1-\beta }\left[ 
\begin{array}{ccc}
3 & \beta z & z \\ 
5\left( 1-\beta \right) +\beta z & 1 & z \\ 
z & z & z%
\end{array}%
\right] .
\end{equation*}%
and the Replicator Dynamics is as follows. 
\begin{equation}
\begin{array}{l}
\dot{x}_{1}=\frac{x_{1}}{1-\beta }\left( \left( z-3\right) x_{1}^{2}+\left(
1-\beta \right) \left( 2z-5\right) x_{1}x_{2}+\left( z-1\right)
x_{2}^{2}+\left( 3-z\right) x_{1}+\left( \beta z-z\right) x_{2}\right)  \\ 
\dot{x}_{2}=\frac{x_{2}}{1-\beta }\left( \left( z-3\right) x_{1}^{2}+\left(
1-\beta \right) \left( 2z-5\right) x_{1}x_{2}+\left( z-1\right)
x_{2}^{2}+\left( 1-\beta \right) \left( 5-z\right) x_{1}+\left( 1-z\right)
x_{2}\right) 
\end{array}
\label{RD2d}
\end{equation}

There are four fixed points for this Dynamics and a standard linearization
analysis produces the results shown in table 2. \medskip

\begin{center}
\begin{tabular}{c}
\begin{tabular}{|c|c|c|}
\hline
Point & eigenvectors & eigenvalues \\ \hline
$\left\{ 0,0\right\} $ & $%
\begin{array}{c}
e_{1}=\left( 1,0\right)  \\ 
\medskip e_{2}=\left( 0,1\right) 
\end{array}%
$ & $%
\begin{array}{c}
\lambda _{1}=0 \\ 
\medskip \lambda _{2}=0%
\end{array}%
$ \\ \hline
$\left\{ 1,0\right\} $ & $%
\begin{array}{c}
e_{1}=\left( 1,0\right)  \\ 
\medskip e_{2}=\left( -1,1\right) 
\end{array}%
$ & $%
\begin{array}{c}
\lambda _{1}=\frac{z-3}{1-\beta } \\ 
\medskip \lambda _{2}=\frac{2-5\beta +\beta z}{1-\beta }%
\end{array}%
$ \\ \hline
$\left\{ 0,1\right\} $ & $%
\begin{array}{c}
e_{1}=\left( 1,-1\right)  \\ 
\medskip e_{2}=\left( 0,1\right) 
\end{array}%
$ & $%
\begin{array}{c}
\lambda _{1}=\frac{\beta z-1}{1-\beta } \\ 
\medskip \lambda _{2}=\frac{z-1}{1-\beta }%
\end{array}%
$ \\ \hline
$\left\{ \frac{\beta z-1}{1+2\beta z-5\beta },\frac{\beta z+2-5\beta }{%
1+2\beta z-5\beta }\right\} $ & $%
\begin{array}{c}
e_{1}=\left( -1,1\right)  \\ 
\medskip e_{2}=\left( 1,\frac{\beta z-5\beta +2}{\beta z-1}\right) 
\end{array}%
$ & $%
\begin{array}{c}
\lambda _{1}=\frac{\left( 1-\beta z\right) \left( \beta z-5\beta +2\right) }{%
\left( 1-\beta \right) \left( 2\beta z-5\beta +1\right) } \\ 
\medskip \lambda _{2}=\frac{z\beta \left( 2-\beta \right) \left( z-5\right)
+3+z}{\left( 1-\beta \right) \left( 2\beta z-5\beta +1\right) }%
\end{array}%
$ \\ \hline
\end{tabular}
\\ 
$\text{\rule{0in}{0.35in}}$Table 2. Eigenvalues and eigenvectors for the
fixed points  \\ 
in the replicator dynamics system given by equations (\ref{RD2d}).$\text{%
\rule[-0.3in]{0in}{0.18in}}$%
\end{tabular}
\end{center}

Depending on the values of the parameters $z$ and $\beta $ we obtain
different solutions for the dynamics~(\ref{RD2d}). The $\beta -z$ parameter
space can be divided into 10 regions (see figure 1) which have qualitatively
different pictures of the dynamics (see figure~2). The main features of this
overall picture can be summarized as follows. If $z>3$ then the population
evolves towards a monomorphic state in which every player uses the
unsociable strategy $\sigma _{B}$. If $z<1$ then the picture resembles the
iterated prisoner's dilemma: if $\beta $ is small then defection is stable,
but if $\beta $ is large populations using either defection or conditional
cooperation are asymptotically stable and the population which arises
depends on the initial conditions. The most interesting dynamics occur when $%
\beta $ is large and $1<z<3$ (labelled as regions VI to IX in figure 1). If $%
\beta $ is large and $z<\frac{5\beta -2}{\beta }$ then conditional
cooperation is the only asymptotically stable behaviour, and in region VII
this is the endpoint of all trajectories which start in the interior of the
simplex. In region VI a polymorphic population is stable: a proportion of
players, $x$, use the conditional cooperative strategy and a proportion, $%
1-x $, defect where $x=\frac{\beta z-1}{1+2\beta z-5\beta }$. In this
population the proportion of individuals that would be observed in
cooperative partnerships is $x^{2}$; the proportion of individuals involved
in partnerships for which mutual defection was the norm would be $(1-x)^{2};$
and a proportion $2x(1-x)$ of individuals would be living alone.

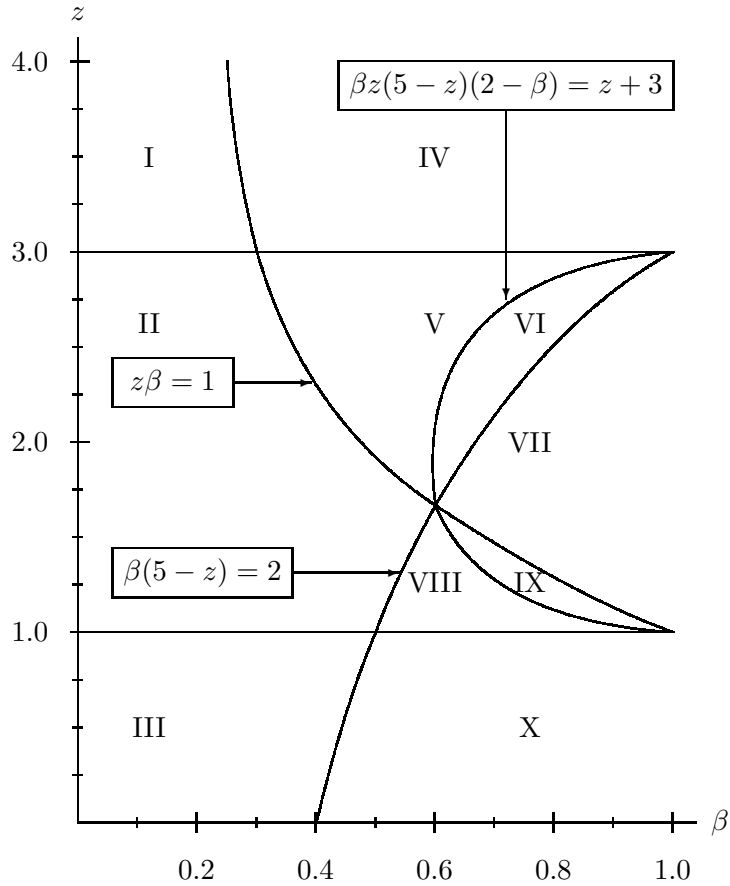
\begin{figure}[ht]
\centering
\begin{minipage}[b]{1.0\textwidth}
\centering
\unitlength=0.90pt
\begin{picture}(265.00,450.00)(0.00,0.00) 
\linethickness{0.2mm} 
\put(300.00,30.00){\makebox(0.00,0.00){$\beta$}} 
\put(30.00,370.00){\makebox(0.00,0.00){$z$}} 
\put(10.00,350.00){\makebox(0.00,0.00){$4.0$}} 
\put(10.00,270.00){\makebox(0.00,0.00){$3.0$}} 
\put(10.00,190.00){\makebox(0.00,0.00){$2.0$}} 
\put(10.00,110.00){\makebox(0.00,0.00){$1.0$}} 
\put(280.00,10.00){\makebox(0.00,0.00){$1.0$}} 
\put(230.00,10.00){\makebox(0.00,0.00){$0.8$}} 
\put(180.00,10.00){\makebox(0.00,0.00){$0.6$}} 
\put(130.00,10.00){\makebox(0.00,0.00){$0.4$}} 
\put(80.00,10.00){\makebox(0.00,0.00){$0.2$}} 
\put(35.00,350.00){\line(-1,0){8.00}} 
\put(35.00,270.00){\line(-1,0){8.00}} 
\put(35.00,190.00){\line(-1,0){8.00}} 
\put(35.00,110.00){\line(-1,0){8.00}} 
\put(32.00,50.00){\line(-1,0){4.00}} 
\put(32.00,70.00){\line(-1,0){4.00}} 
\put(32.00,90.00){\line(-1,0){4.00}} 
\put(32.00,130.00){\line(-1,0){4.00}} 
\put(32.00,150.00){\line(-1,0){4.00}} 
\put(32.00,170.00){\line(-1,0){4.00}} 
\put(32.00,210.00){\line(-1,0){4.00}} 
\put(32.00,230.00){\line(-1,0){4.00}} 
\put(32.00,250.00){\line(-1,0){4.00}} 
\put(32.00,290.00){\line(-1,0){4.00}} 
\put(32.00,310.00){\line(-1,0){4.00}} 
\put(32.00,330.00){\line(-1,0){4.00}} 
\put(280.00,34.00){\line(0,-1){8.00}} 
\put(230.00,34.00){\line(0,-1){8.00}} 
\put(180.00,34.00){\line(0,-1){8.00}} 
\put(130.00,34.00){\line(0,-1){8.00}} 
\put(80.00,34.00){\line(0,-1){8.00}} 
\put(255.00,32.00){\line(0,-1){4.00}} 
\put(205.00,32.00){\line(0,-1){4.00}} 
\put(155.00,32.00){\line(0,-1){4.00}} 
\put(105.00,32.00){\line(0,-1){4.00}} 
\put(55.00,32.00){\line(0,-1){4.00}} 
\put(30.00,30.00){\line(1,0){260.00}} 
\put(30.00,30.00){\line(0,1){330.00}} 
\put(30.00,110.00){\line(1,0){250.00}} 
\put(30.00,270.00){\line(1,0){250.00}} 
\qbezier(92.50,350.00)(95,310)(105.00,270.00) 
\qbezier(105.00,270.00)(125.00,200.00)(180.00,163.33) 
\qbezier(180.00,163.33)(238.00,125.00)(280.00,110.00) 
\qbezier(280.00,270.00)(225,240)(180.00,163.33) 
\qbezier(180.00,163.33)(162,130)(155.00,110.00) 
\qbezier(155.00,110.00)(142,80)(130.00,30.00) 
\qbezier(280.00,270.00)(168,260)(180.00,163.33) 
\qbezier(180.00,163.33)(200.00,115.00)(280.00,110.00) 
\put(60.00,310.00){\makebox(0.00,0.00){I}} 
\put(60.00,240.00){\makebox(0.00,0.00){II}} 
\put(60.00,70.00){\makebox(0.00,0.00){III}} 
\put(180.00,310.00){\makebox(0.00,0.00){IV}} 
\put(180.00,240.00){\makebox(0.00,0.00){V}} 
\put(220.00,240.00){\makebox(0.00,0.00){VI}} 
\put(220.00,190.00){\makebox(0.00,0.00){VII}} 
\put(180.00,131.00){\makebox(0.00,0.00){VIII}} 
\put(220.00,131.00){\makebox(0.00,0.00){IX}} 
\put(220.00,70.00){\makebox(0.00,0.00){X}} 
\put(45.00,205.00){\framebox(50.00,20.00){$z \beta=1$}} 
\put(95.00,215.00){\vector(1,0){33.00}} 
\put(45.00,125.00){\framebox(75.00,20.00){$ \beta (5-z)=2$}} 
\put(120.00,135.00){\vector(1,0){45.00}} 
\put(140.00,330.00){\framebox(140.00,20.00){$ \beta z (5-z) (2-\beta)=z+3$}} 
\put(210.00,330.00){\vector(0,-1){80.00}} 
\end{picture}
\end{minipage}
\caption{Regions for the parameters $z$ and $\beta$ which lead to qualitatively 
different pictures of the evolutionary dynamics.}
\end{figure}

\begin{figure}[ht]
\centering
\begin{minipage}[b]{0.3\textwidth}
\centering
\unitlength=0.37pt
\begin{picture}(260.00,260.00)(0.00,0.00) 
\thicklines 
\put(230.00,10.00){\makebox(0.00,0.00){$1$}} 
\put(140.00,10.00){\makebox(0.00,0.00){$x_1$}} 
\put(40.00,10.00){\makebox(0.00,0.00){$0$}} 
\put(10.00,230.00){\makebox(0.00,0.00){$1$}} 
\put(10.00,140.00){\makebox(0.00,0.00){$x_2$}} 
\put(10.00,40.00){\makebox(0.00,0.00){$0$}} 
\put(80.00,145.00){\makebox(0.00,0.00){d}} 
\put(180.00,50.00){\makebox(0.00,0.00){b}} 
\put(50.00,180.00){\makebox(0.00,0.00){c}} 
\put(50.00,50.00){\makebox(0.00,0.00){a}} 
\put(35.00,30.00){\line(1,0){190.00}} 
\put(30.00,35.00){\line(0,1){190.00}} 
\put(35.00,225.00){\line(1,-1){44.00}} 
\put(225.00,35.00){\line(-1,1){137.00}} 
\put(83.00,177.00){\circle{10.00}} 
\put(30.00,30.00){\circle{10.00}} 
\put(30.00,230.00){\circle{10.00}} 
\put(230.00,30.00){\circle{10.00}} 
\end{picture}
\end{minipage}%
\hspace{0.04\textwidth}%
\begin{minipage}[b]{0.3\textwidth}
\centering
\subfigure{\includegraphics[totalheight=0.16\textheight, width=3.5cm]{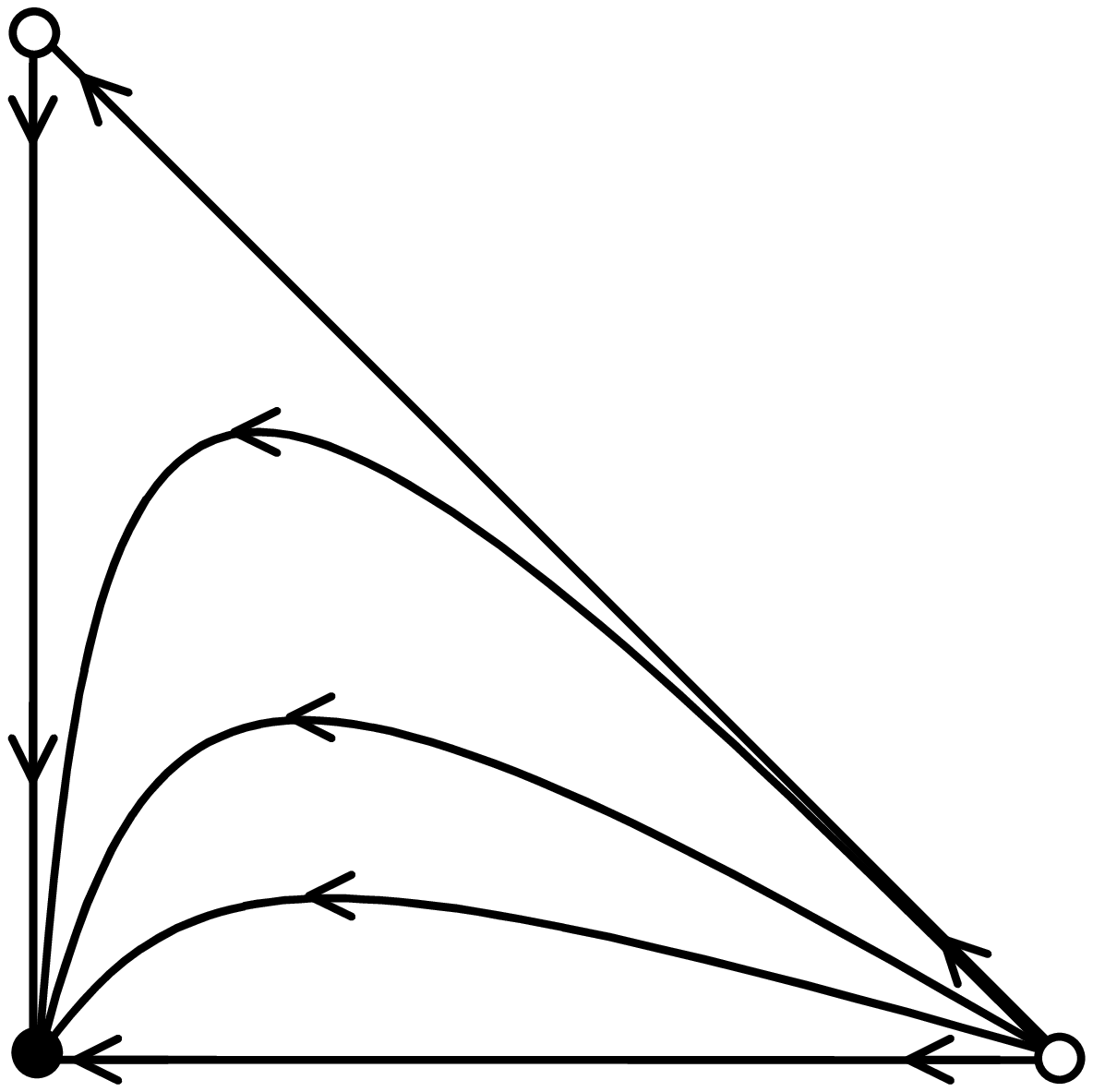}}
I
\end{minipage}%
\hspace{0.03\textwidth}%
\begin{minipage}[b]{0.3\textwidth}
\centering
\subfigure{\includegraphics[totalheight=0.16\textheight, width=3.5cm]{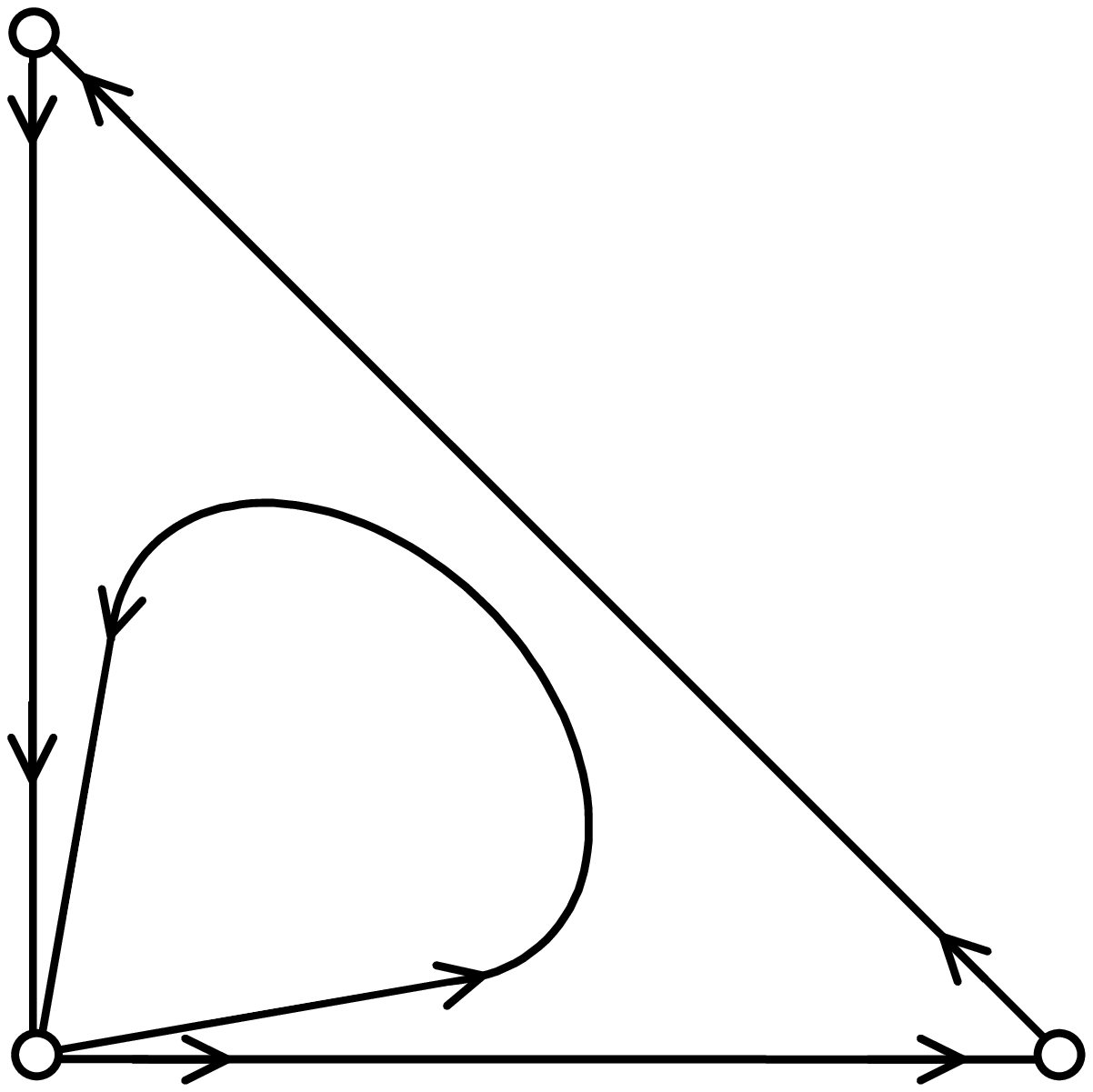}}
II
\end{minipage}\\
\vspace{-0.5cm}
\begin{minipage}[b]{0.3\textwidth}
\centering
\subfigure{\includegraphics[totalheight=0.16\textheight, width=3.5cm]{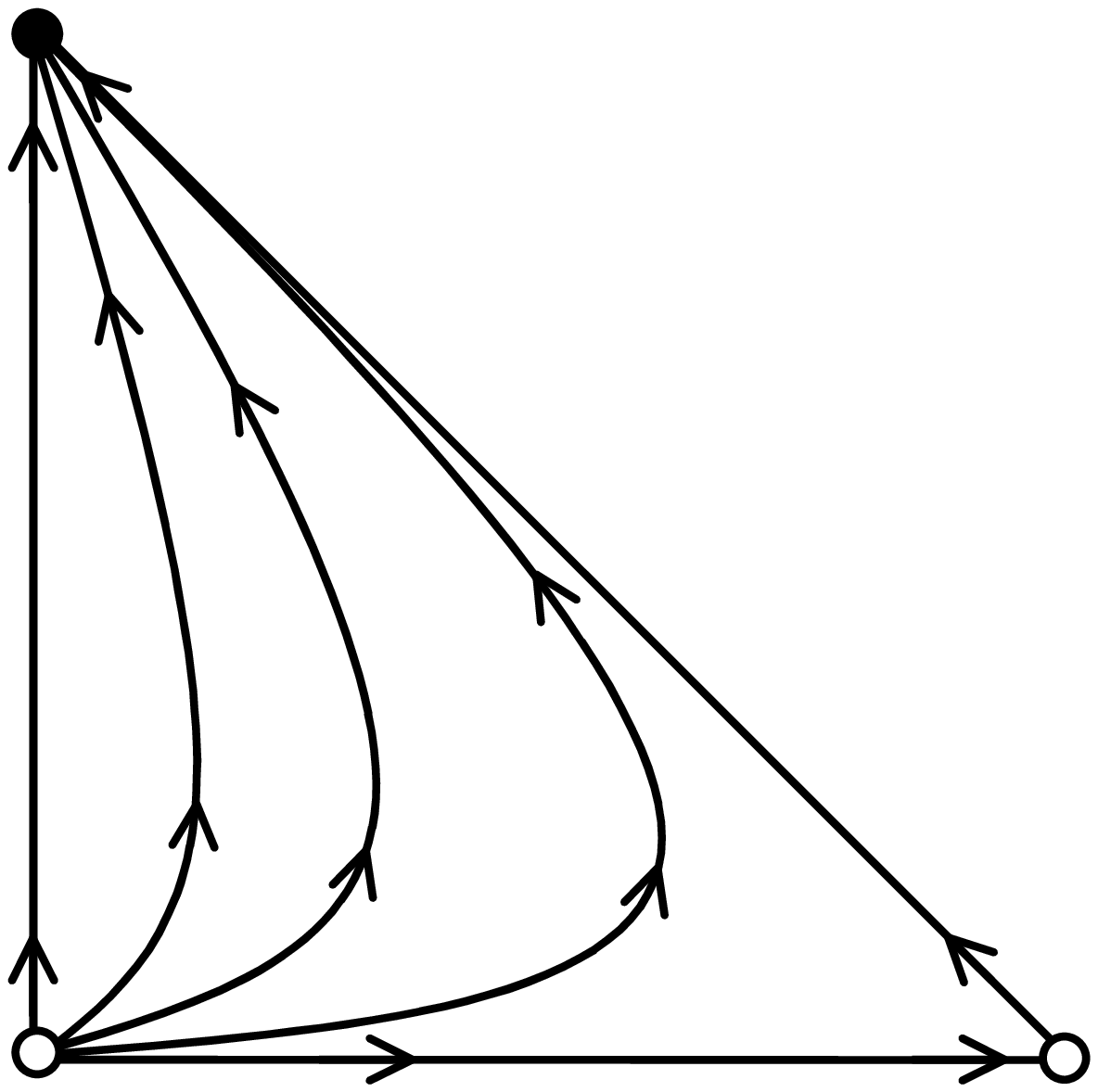}}
III
\end{minipage}%
\hspace{0.04\textwidth}%
\begin{minipage}[b]{0.3\textwidth}
\centering
\subfigure{\includegraphics[totalheight=0.16\textheight, width=3.5cm]{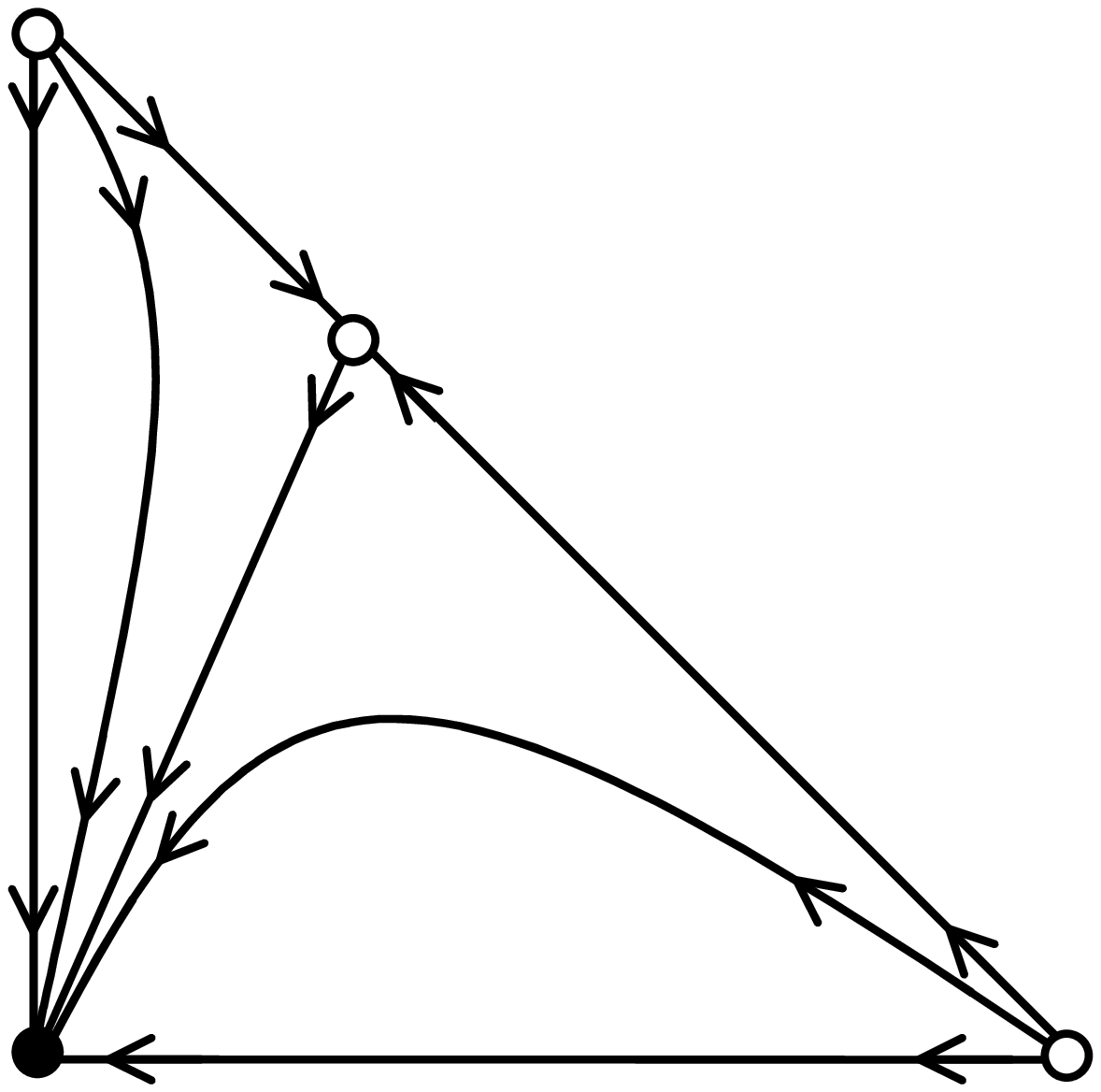}}
IV
\end{minipage}%
\hspace{0.04\textwidth}%
\begin{minipage}[b]{0.3\textwidth}
\centering
\subfigure{\includegraphics[totalheight=0.16\textheight, width=3.5cm]{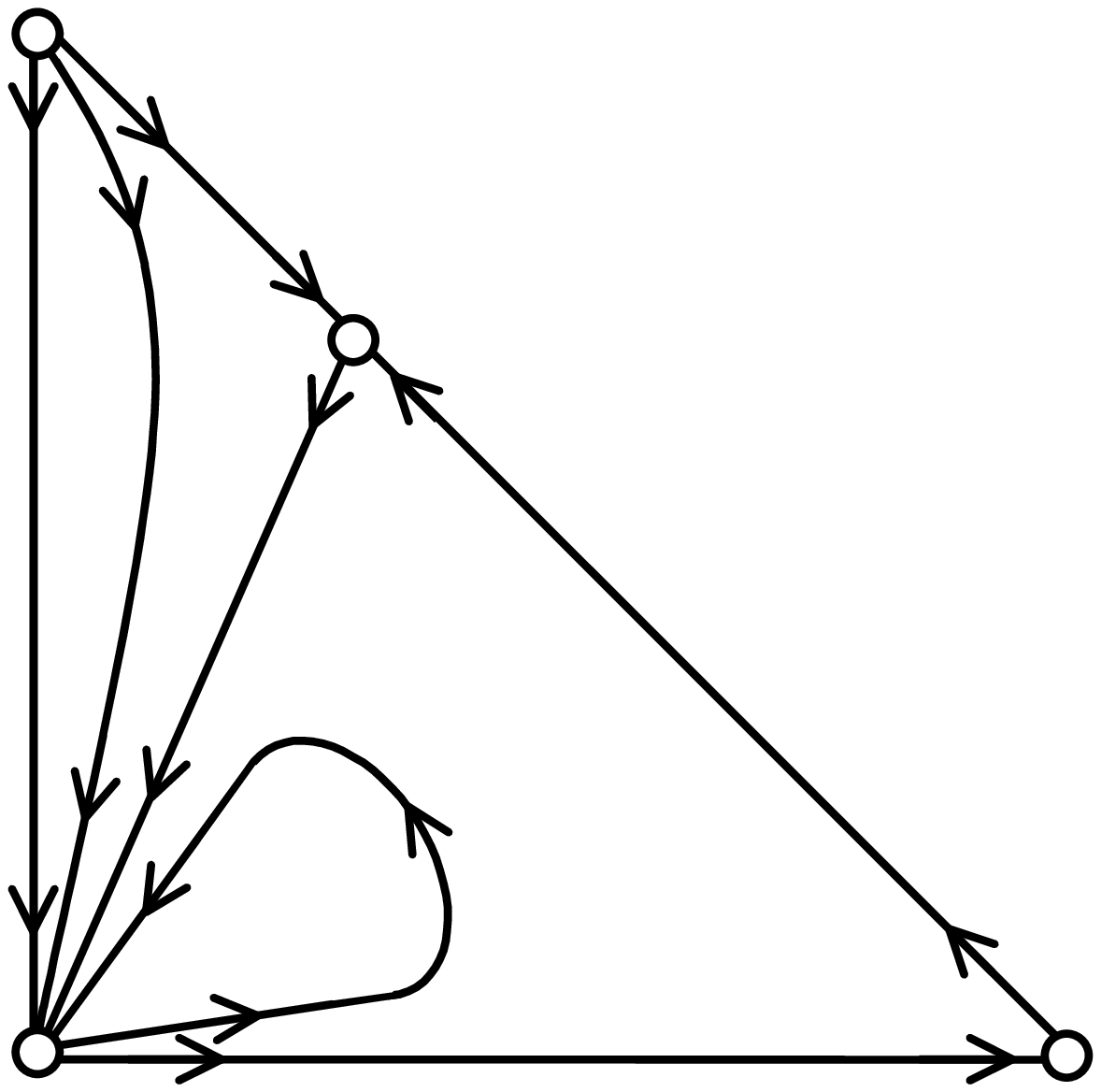}}
V
\end{minipage}\\
\vspace{-0.5cm}
\begin{minipage}[b]{0.3\textwidth}
\centering
\subfigure{\includegraphics[totalheight=0.16\textheight, width=3.5cm]{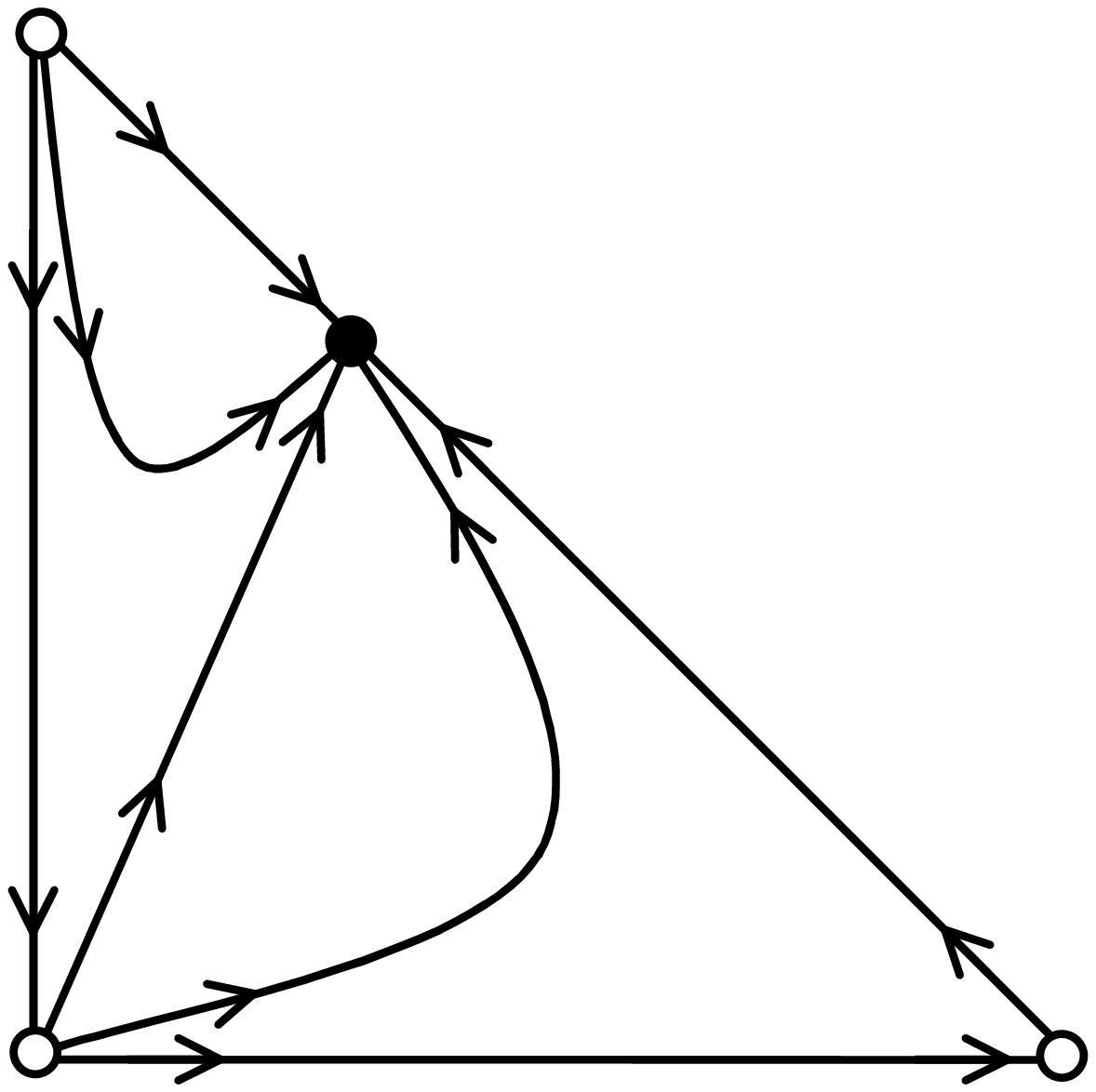}}
VI
\end{minipage}%
\hspace{0.04\textwidth}%
\begin{minipage}[b]{0.3\textwidth}
\centering
\subfigure{\includegraphics[totalheight=0.16\textheight, width=3.5cm]{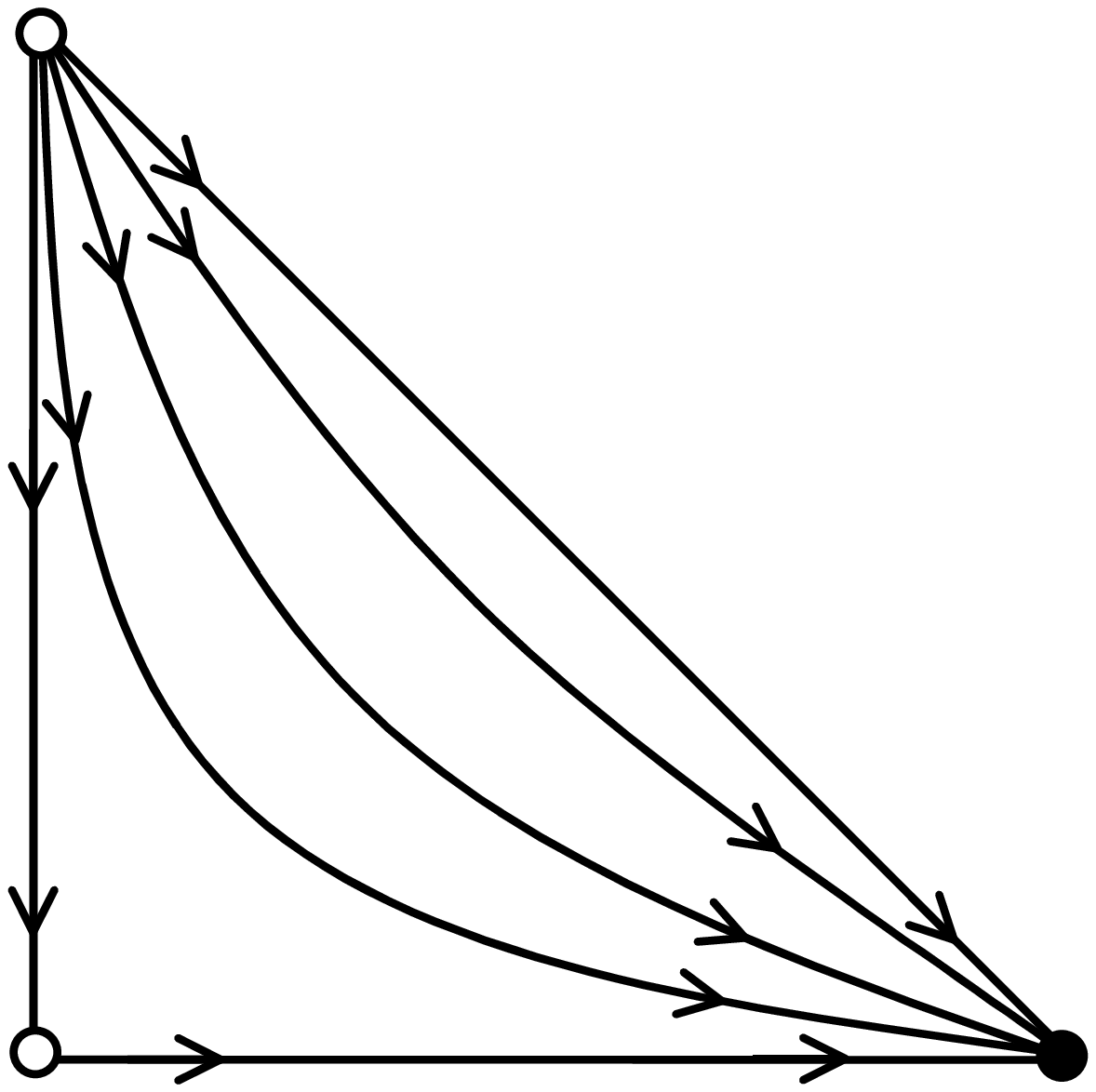}}
VII
\end{minipage}%
\hspace{0.04\textwidth}%
\begin{minipage}[b]{0.3\textwidth}
\centering
\subfigure{\includegraphics[totalheight=0.16\textheight, width=3.5cm]{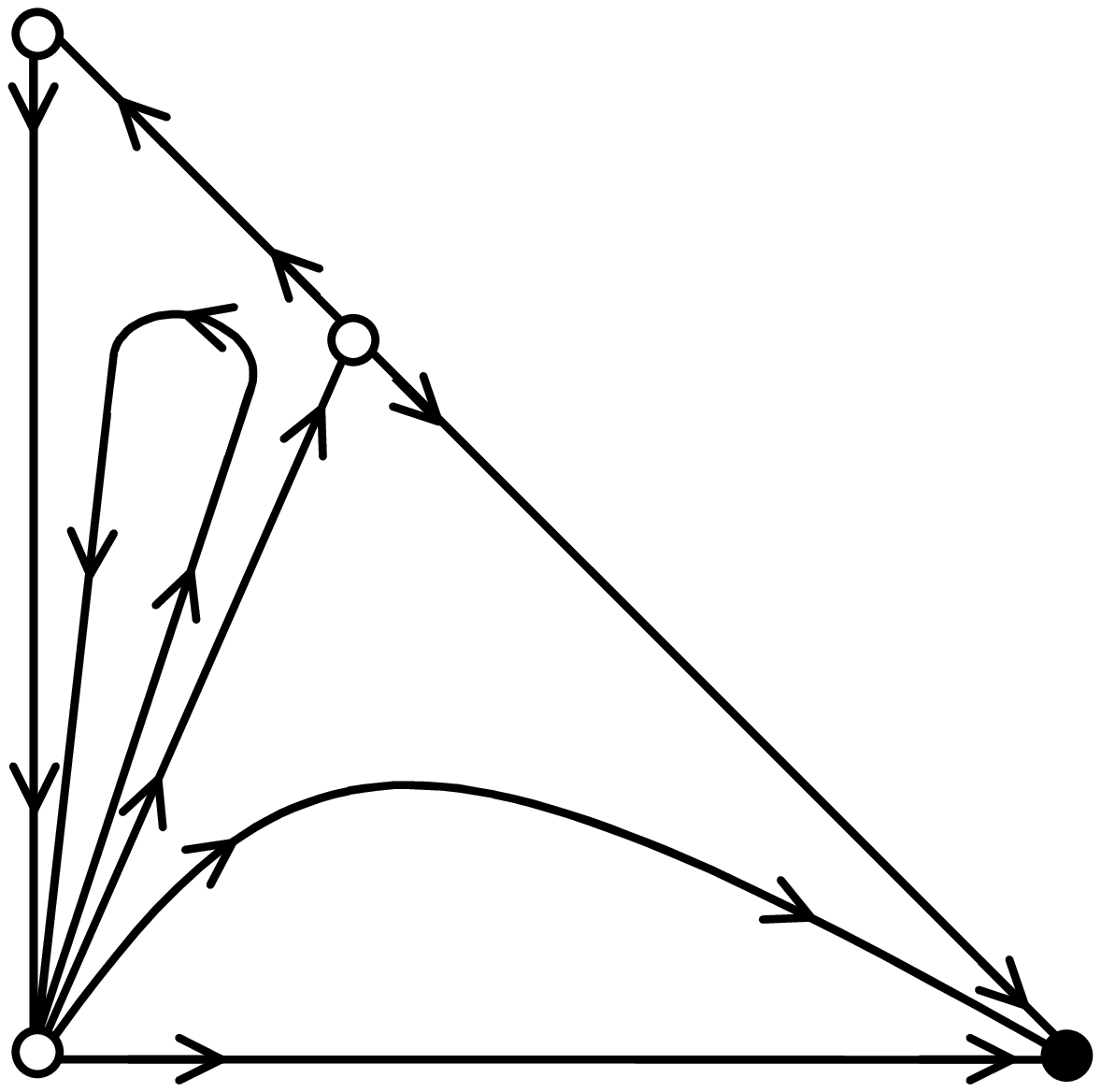}}
VIII
\end{minipage}\\
\vspace{-0.5cm}
\begin{minipage}[b]{0.3\textwidth}
\centering
\subfigure{\includegraphics[totalheight=0.16\textheight, width=3.5cm]{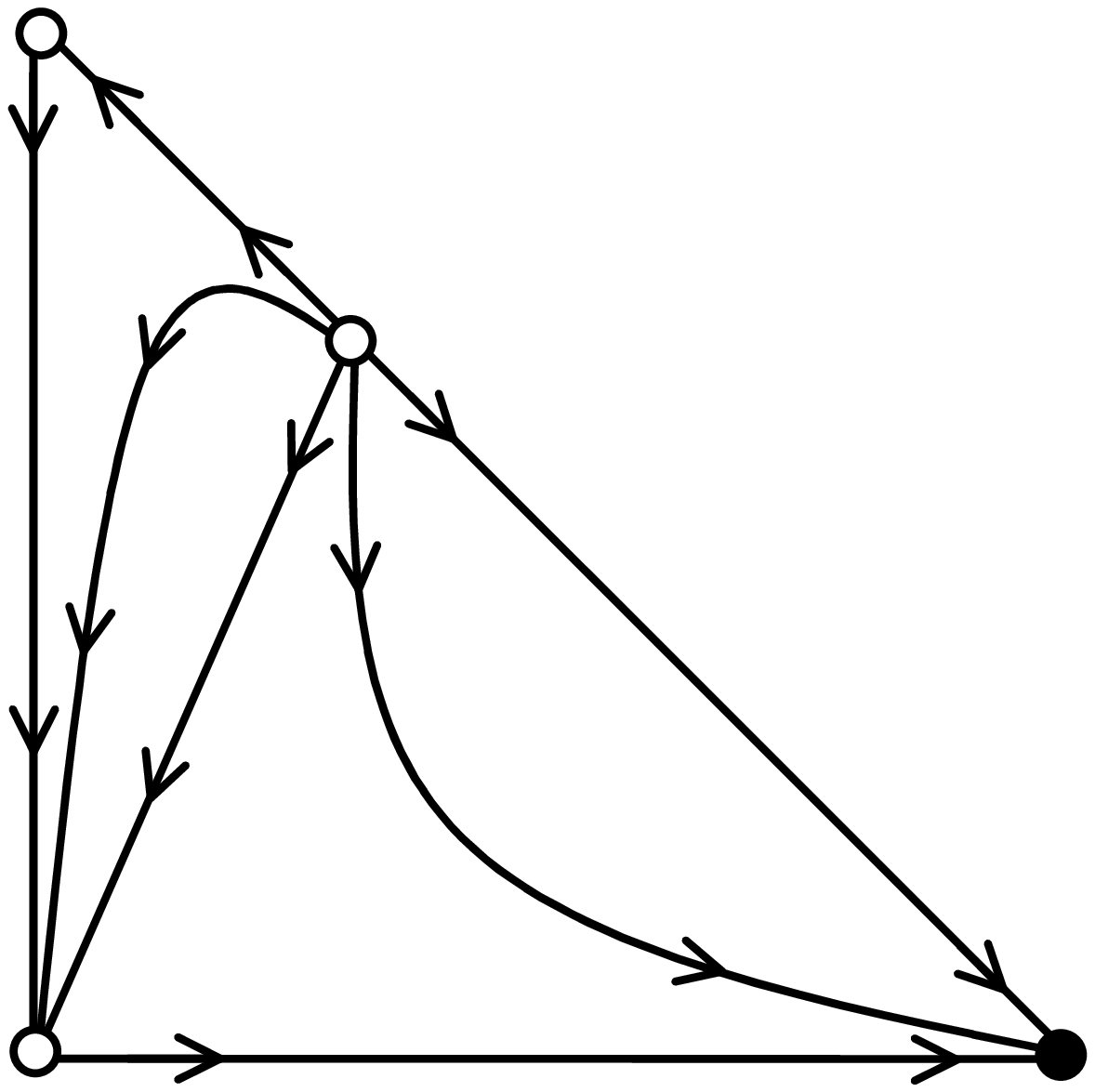}}
IX
\end{minipage}%
\hspace{0.04\textwidth}%
\begin{minipage}[b]{0.3\textwidth}
\centering
\subfigure{\includegraphics[totalheight=0.16\textheight, width=3.5cm]{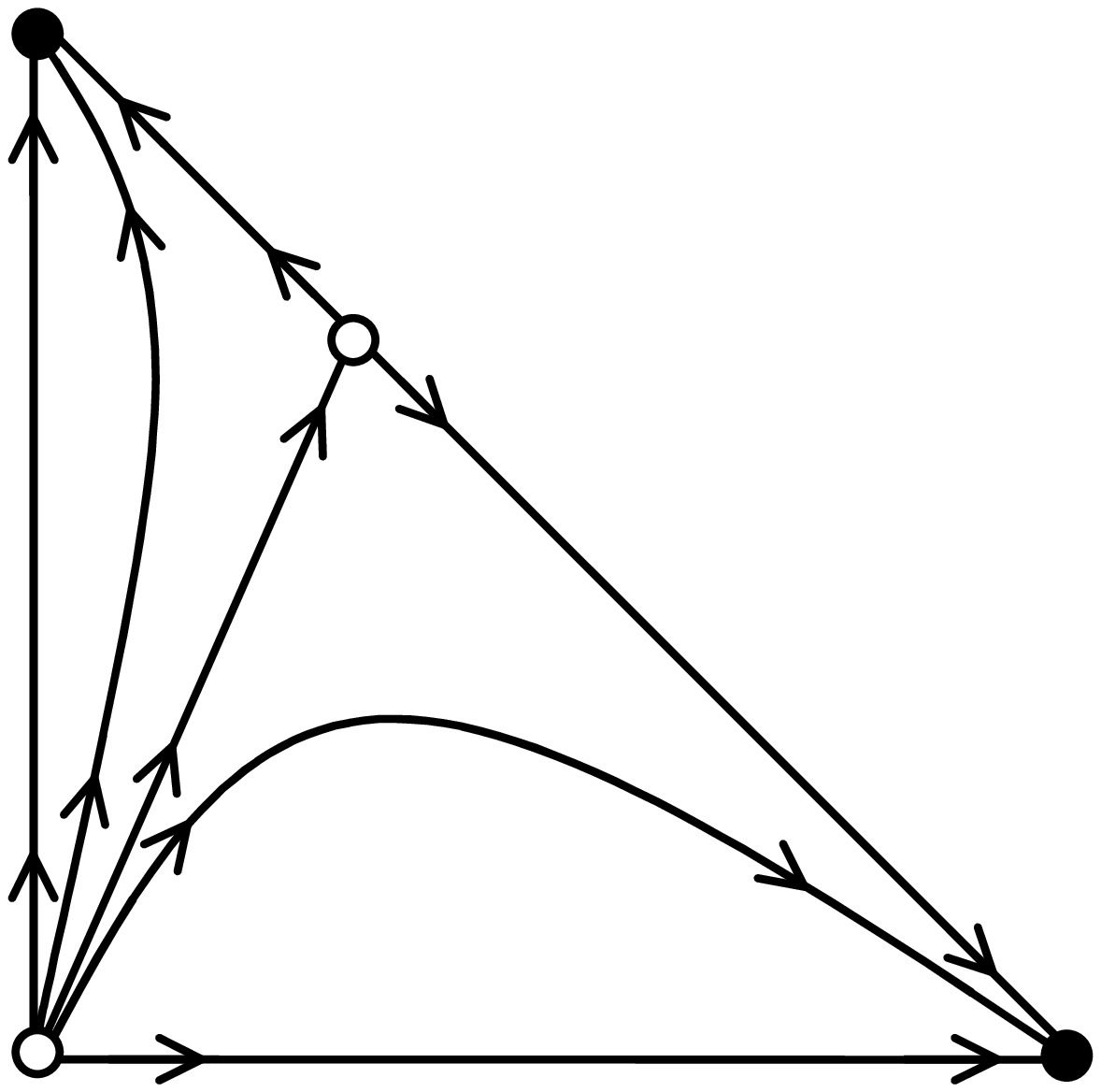}}
X
\end{minipage}
\caption{Qualitative pictures of the dynamics in the different regions of
the $\protect\beta-z$ parameter space. The regions are labelled according to
figure 1. Point $a$ corresponds to a population which consists of 100\% of
players using strategy $\protect\sigma_{B}$. Point $b$ corresponds to a
population which consists of 100\% of players using strategy $\protect\sigma%
_{C}$. Point $c$ corresponds to a population which consists of 100\% of
players using strategy $\protect\sigma_{D}$. Point $d$ corresponds to a
polymorphic population with players using either $\protect\sigma_{C}$ or $%
\protect\sigma_{D}$. Asymptotically stable points are shown as solid
circles, the other fixed points are shown as open circles.}
\end{figure}

\section*{Introducing unconditional cooperators}

In region VII of the $\beta -z$ parameter space we have found that
conditional cooperation is asymptotically stable. It is pertinent to ask
whether this property would be destroyed if we allowed individuals to use
the ``sucker'' strategy of unconditional cooperation (which we denote $%
\sigma _{S}$). Recall that in the iterated prisoner's dilemma, tit-for-tat
is not asymptotically stable due to the presence of unconditional
cooperators. Similarly, it is conceivable that the polymorphic population
which is stable in region VI could be destabilized by the introduction of a
strategy of unconditional cooperation.

We introduce a proportion $x_{3}$ of players who use the strategy of
unconditional cooperation, $\sigma _{S}$. These players associate in $G_{0}$
and cooperate in $G_{1}$ whatever their opponent does. (The proportion of
individuals using $\sigma _{B}$ is then $1-x_{1}-x_{2}-x_{3}$.) The new
payoff matrix is given by%
\begin{equation}
A=\left[ 
\begin{array}{llll}
\pi (\sigma _{C},\sigma _{C}) & \pi (\sigma _{C},\sigma _{D}) & \pi (\sigma
_{C},\sigma _{S}) & \pi (\sigma _{C},\sigma _{B}) \\ 
\pi (\sigma _{D},\sigma _{C}) & \pi (\sigma _{D},\sigma _{D}) & \pi (\sigma
_{D},\sigma _{S}) & \pi (\sigma _{D},\sigma _{B}) \\ 
\pi (\sigma _{S},\sigma _{C}) & \pi (\sigma _{S},\sigma _{D}) & \pi (\sigma
_{S},\sigma _{S}) & \pi (\sigma _{S},\sigma _{B}) \\ 
\pi (\sigma _{B},\sigma _{C}) & \pi (\sigma _{B},\sigma _{D}) & \pi (\sigma
_{B},\sigma _{S}) & \pi (\sigma _{B},\sigma _{B})%
\end{array}%
\right] =\frac{1}{1-\beta }\left[ 
\begin{array}{cccc}
3 & \beta z & 3 & z \\ 
5\left( 1-\beta \right) +\beta z & 1 & 5 & z \\ 
3 & 0 & 3 & z \\ 
z & z & z & z%
\end{array}%
\right]  \label{RD3d}
\end{equation}

An analysis of the corresponding Replicator dynamics leads to the dynamics
shown in figures 3 and 4 for regions VI and VII respectively (see appendix B
for details). These figures show the dynamics for particular values of $z$
and $\beta $ but the pictures are qualitatively similar for any values of
these parameters in the appropriate range. From these figures we can see
that the polymorphic population remains asymptotically stable in region VI.
In region VII, the population of conditional cooperators is no longer
asymptotically stable. However, populations which consist of mixtures of
conditional and unconditional cooperators are the only end points of all
solution trajectories which start in the interior of the simplex.

\begin{figure}[ht]
\centering
\begin{minipage}[b]{0.45\textwidth}
\centering
\subfigure{\includegraphics[totalheight=0.3\textheight, width=7cm]{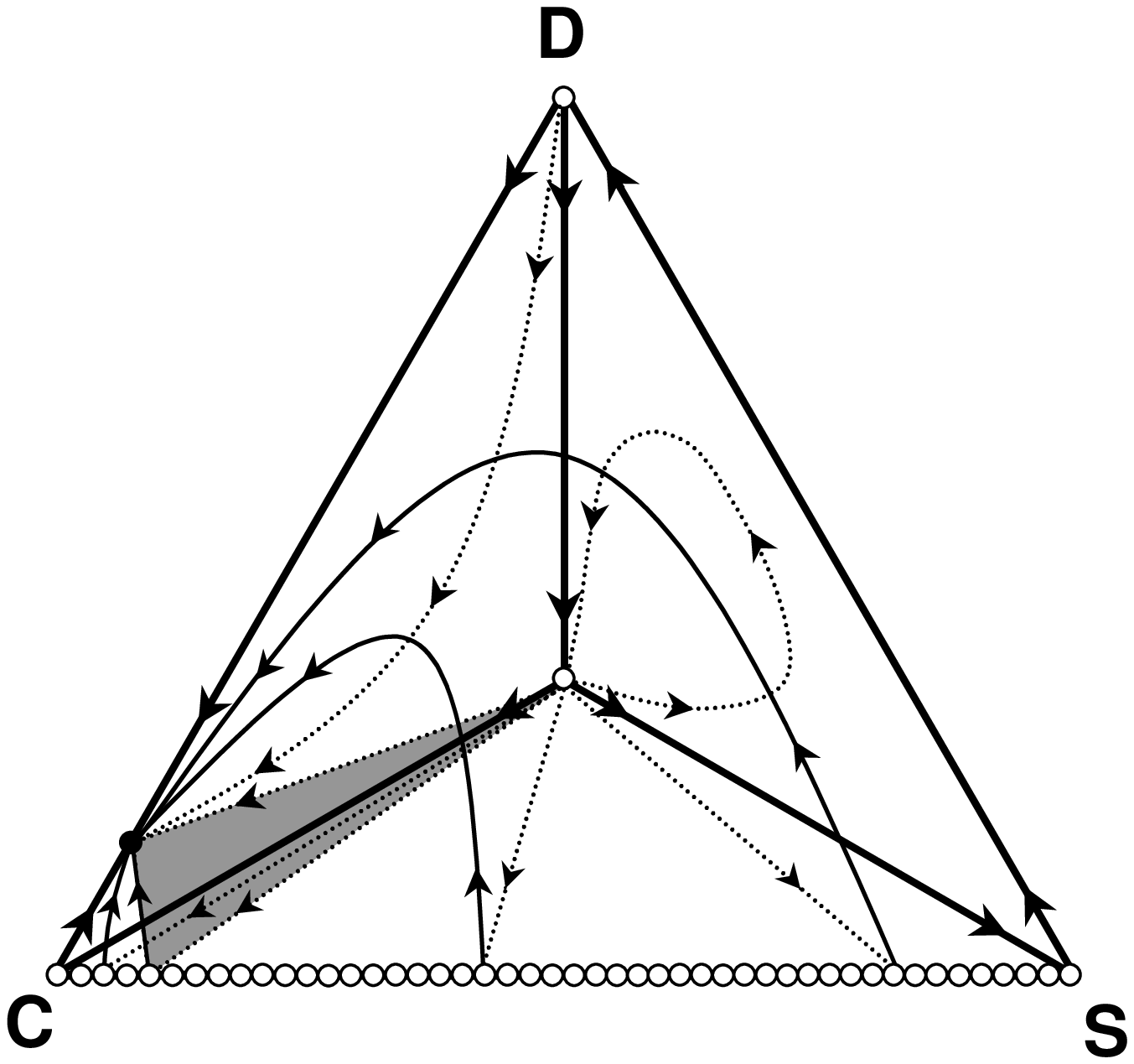}}
Region VI: $z=\frac{25}{10}$ and $\protect\beta = \frac{75}{100}$. 
\end{minipage}%
\hspace{0.04\textwidth}%
\begin{minipage}[b]{0.45\textwidth}
\centering
\subfigure{\includegraphics[totalheight=0.3\textheight, width=7cm]{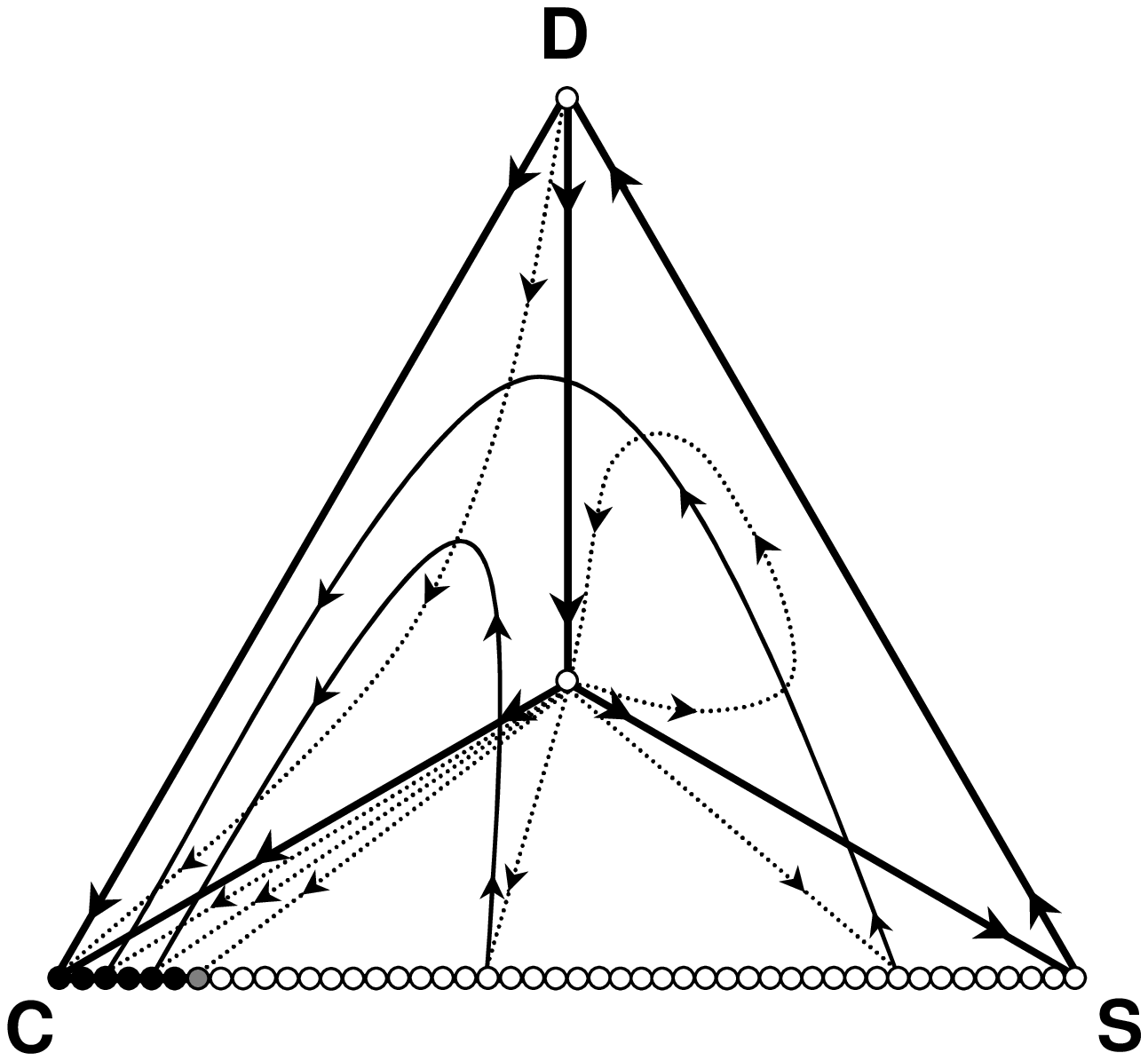}}
Region VII: $z=\frac{25}{10}$ and $\protect\beta =\frac{9}{10}$.
 \end{minipage}
\caption{Qualitative picture of the dynamics for the replicator system with
payoff matrix given by equation (\ref{RD3d}). Vertex $S$ corresponds to a
population which consists of 100\% of players using strategy $\protect\sigma%
_{S}$. Vertex $C$ corresponds to a population which consists of 100\% of
players using strategy $\protect\sigma_{C}$. Vertex $D$ corresponds to a
population which consists of 100\% of players using strategy $\protect\sigma%
_{D}$. The unlabelled vertex corresponds to a population in which all
players live alone.}
\end{figure}

\section*{Discussion}

A minimal version of the iterated prisoner's dilemma deals with a population
consisting of unconditional cooperators, unconditional defectors and
conditional cooperators (such as tit-for-tat). In that model there is a
threshold problem: cooperative behaviour only evolves if the initial
proportion of conditional cooperators exceeds some value \cite{VR, D}.
Although it is sometimes suggested that the always defect strategy is an ESS
or that the corresponding population is asymptotically stable, this is not
the case. If sufficiently many varied strategies are introduced then the
barrier can be removed \cite{NS}.

We have introduced an outside option into the iterated prisoner's dilemma,
which allows individuals to avoid being condemned to maintain an
unprofitable interaction of permanent mutual defection. This provides
another way of removing the barrier to the evolution of cooperative
behaviour. The requirement is that the payoff from the outside option should
be neither so poor that it is irrelevant nor so high that everyone opts for
a solitary existence. The existence of the outside option also admits a
range of parameter values for which a polymorphic population involving
defectors and conditional cooperators is asymptotically stable, even in the
presence of unconditional cooperators.

Some of the results we have obtained are similar to those obtained for
optional public good games, which are multi-player generalizations of the
prisoner's dilemma \cite{HMHS}. In these games, as in ours, making
participation voluntary enhances the possibilities for cooperation. One
difference between the two models is that in the optional public good game
rock-scissors-paper style cycles may occur. In our model, such cyclic
behaviour does not arise. However, in both models the fixed point
representing non-participatory behaviour may be non-hyperbolic. This leads
to periods of cooperative behaviour, but eventually the population returns
to a state in which everyone lives alone.

The iterated prisoner's dilemma is an unrealistic model of social
interactions because it treats one type of interaction between individuals
in isolation from all others. We have shown, by means of a relatively simple
example, that the methods of stochastic game theory can be employed to
overcome this restriction. The prisoners dilemma has also been criticized as
being an unrealistic model of social interactions on other grounds \cite{DMG}%
. Our approach is not specific to the prisoner's dilemma. That context game
may be replaced by any other game or, indeed, a game which is randomly
selected with a known probability from a set of games \cite{K}. This allows
quite complex social behaviour to be analyzed.

\section*{Appendix A}

To integrate the Replicator Dynamics system (\ref{RD2d}) we make the
following coordinate substitutions. 
\begin{equation*}
k=\frac{x_{2}}{x_{1}}\ \ \text{\ and\ }\ \ l=\frac{1-x_{1}-x_{2}}{x_{1}}
\end{equation*}
or 
\begin{equation*}
x_{1}=\frac{1}{1+k+l}\text{ \ \ \ \ \ and \ \ \ \ \ }x_{2}=\frac{k}{1+k+l}.
\end{equation*}
Then we have 
\begin{equation*}
\begin{array}{l}
\dot{k}=k\frac{a+bk}{1+l+k} \\ 
\dot{l}=l\frac{c+fk}{1+l+k}%
\end{array}
.
\end{equation*}
Solution trajectories can be found by integrating 
\begin{equation*}
\frac{dl}{dk}=\frac{l\left( c+fk\right) }{k\left( a+bk\right) }
\end{equation*}

\begin{equation*}
\frac{dl}{l}=\frac{c}{a}\frac{dk}{k}+\frac{af-bc}{a}\frac{dk}{\left(
a+bk\right) }
\end{equation*}

This can be done analytically to obtain 
\begin{eqnarray*}
bc\ln \left| k\right| +\left( af-bc\right) \ln \left| a+bk\right| &=&ab\ln
\left| l\right| +C \\
\left| k\right| ^{bc}\left| a+bk\right| ^{\left( af-bc\right) } &=&C\left|
l\right| ^{ab}
\end{eqnarray*}
where $C$ is\ a constant that depends on the initial conditions. Finally,
substituting the expressions for $k$ and $l$ into the above formula, we find
that the solution trajectories are described by the expression. 
\begin{equation*}
\left| \frac{x_{2}}{x_{1}}\right| ^{bc}\left| a+b\frac{x_{2}}{x_{1}}\right|
^{\left( af-bc\right) }=C\left| \frac{1-x_{1}-x_{2}}{x_{1}}\right| ^{ab}.
\end{equation*}

\section*{Appendix B}

The Replicator Dynamics with payoff matrix (\ref{RD3d}) is given by the
following system of equations. 
\begin{equation*}
\begin{array}{l}
\dot{x}_{1}=\frac{x_{1}\left( \left( x_{1}+x_{3}\right) \left(
1-x_{1}-x_{3}\right) \left( 3-z\right) +\left( 2z-5\right) \left( \left(
1-\beta \right) x_{1}+x_{3}\right) x_{2}+\left( \left( z-1\right)
x_{2}-\left( 1-\beta \right) z\right) x_{2}\right) }{1-\beta } \\ 
\dot{x}_{2}=\frac{x_{2}\left( \left( x_{1}+x_{3}\right) ^{2}\left(
z-3\right) +\left( \left( 5-z\right) +\left( 2z-5\right) x_{2}\right) \left(
\left( 1-\beta \right) x_{1}+x_{3}\right) +\left( 1-z\right) \left(
1-x_{2}\right) x_{2}\right) }{1-\beta } \\ 
\dot{x}_{3}=\frac{x_{3}\left( \left( x_{1}+x_{3}\right) \left(
1-x_{1}-x_{3}\right) \left( 3-z\right) +\left( 2z-5\right) \left( \left(
1-\beta \right) x_{1}+x_{3}\right) x_{2}+\left( \left( z-1\right)
x_{2}-z\right) x_{2}\right) }{1-\beta }%
\end{array}%
\end{equation*}
The fixed points together with their associated eigenvectors and eigenvalues
are given in table 3. \medskip

\begin{center}
\begin{tabular}{c}
\begin{tabular}{|c|c|c|c|}
\hline
Population & Point & Eigenvectors & Eigenvalues \\ \hline
$100\%$ $\text{of}\ \sigma _{B}$ & $_{\left\{ 0,0,0\right\} }$ & $%
\begin{array}{c}
_{e_{1}=\left( 1,0,0\right) } \\ 
_{e_{2}=\left( 0,1,0\right) } \\ 
_{\medskip e_{2}=\left( 0,0,1\right) }%
\end{array}%
$ & $%
\begin{array}{c}
_{\lambda _{1}=0} \\ 
_{\lambda _{2}=0} \\ 
\medskip _{\lambda _{3}=0}%
\end{array}%
$ \\ \hline
$100\%\text{ of}\ \sigma _{D}$ & $_{\left\{ 0,1,0\right\} }$ & $%
\begin{array}{c}
_{e_{1}=\left( -1,1,0\right) } \\ 
_{e_{2}=\left( 0,1,0\right) } \\ 
\medskip _{e_{2}=\left( 0,1,-1\right) }%
\end{array}%
$ & $%
\begin{array}{c}
_{\lambda _{1}=\frac{\beta z-1}{1-\beta }} \\ 
_{\lambda _{2}=\frac{z-1}{1-\beta }} \\ 
\medskip _{\lambda _{3}=-\frac{1}{1-\beta }}%
\end{array}%
$ \\ \hline
$100\%\text{ of}\ \sigma _{C}$ & $_{\left\{ 1,0,0\right\} }$ & $%
\begin{array}{c}
_{e_{1}=\left( 1,0,-1\right) } \\ 
_{e_{2}=\left( -1,1,0\right) } \\ 
\medskip _{e_{3}=\left( 1,0,0\right) }%
\end{array}%
$ & $%
\begin{array}{c}
_{\lambda _{1}=0} \\ 
_{\lambda _{2}=\frac{2-5\beta +\beta z}{1-\beta }} \\ 
\medskip _{\lambda _{3}=\frac{z-3}{1-\beta }}%
\end{array}%
$ \\ \hline
$100\%\text{ of}\ \sigma _{S}$ & $_{\left\{ 0,0,1\right\} }$ & $%
\begin{array}{c}
_{e_{1}=\left( 1,0,-1\right) } \\ 
_{e_{2}=\left( 0,1,-1\right) } \\ 
\medskip _{e_{3}=\left( 0,0,1\right) }%
\end{array}%
$ & $%
\begin{array}{c}
_{\lambda _{1}=0} \\ 
_{\lambda _{2}=\frac{2}{1-\beta }} \\ 
\medskip _{\lambda _{3}=\frac{z-3}{1-\beta }}%
\end{array}%
$ \\ \hline
$%
\begin{array}{c}
\alpha \%\text{ of }\sigma _{C} \\ 
\text{and} \\ 
\left( 1-\alpha \right) \%\text{ of}\ \sigma _{S}%
\end{array}%
$ & $_{\left\{ \alpha ,0,1-\alpha \right\} }$ & $%
\begin{array}{c}
_{e_{1}=\left( 1,0,-1\right) } \\ 
_{e_{2}=\left( \alpha \frac{\beta z-2\alpha \beta z+5\alpha \beta -2}{%
2-5\alpha \beta +\alpha \beta z},1,\frac{\left( \alpha -1\right) \left(
2\alpha \beta z-5\alpha \beta +2\right) }{2-5\alpha \beta +\alpha \beta z}%
\right) } \\ 
\medskip _{e_{3}=\left( 1,0,\frac{1-\alpha }{\alpha }\right) }%
\end{array}%
$ & $%
\begin{array}{c}
_{\lambda _{1}=0} \\ 
_{\lambda _{2}=\frac{2-5\alpha \beta +\alpha \beta z}{1-\beta }} \\ 
\medskip _{\lambda _{3}=\frac{z-3}{1-\beta }}%
\end{array}%
$ \\ \hline
$%
\begin{array}{c}
\frac{\beta z-1}{2\beta z-5\beta +1}\%\text{ of }\sigma _{C} \\ 
\text{and} \\ 
\frac{2-5\beta +\beta z}{2\beta z-5\beta +1}\%\text{ of}\ \sigma _{S}%
\end{array}%
$ & $_{\left\{ \frac{\beta z-1}{2\beta z-5\beta +1},\frac{2-5\beta +\beta z}{%
2\beta z-5\beta +1},0\right\} }$ & $%
\begin{array}{c}
_{e_{1}=\left( -1,1,0\right) } \\ 
_{e_{2}=\left( 1,\frac{2-5\beta +\beta z}{\beta z-1},0\right) } \\ 
\medskip _{e_{3}=\left( 1,\beta \frac{5-3z}{\beta z-1},\frac{2\beta z-5\beta
+1}{\beta z-1}\right) }%
\end{array}%
$ & $%
\begin{array}{c}
_{\lambda _{1}=\frac{\left( 1-\beta z\right) \left( \beta z-5\beta +2\right) 
}{\left( 1-\beta \right) \left( 2\beta z-5\beta +1\right) }} \\ 
_{\lambda _{2}=\frac{z\beta \left( 2-\beta \right) \left( z-5\right) +3+z}{%
\left( 1-\beta \right) \left( 2\beta z-5\beta +1\right) }} \\ 
\medskip _{\lambda _{3}=\frac{-\beta z\left( 2-5\beta +\beta z\right) }{%
\left( 1-\beta \right) \left( 2\beta z-5\beta +1\right) }}%
\end{array}%
$ \\ \hline
\end{tabular}
\\ 
$\text{\rule{0in}{0.35in}}$Table 3. Eigenvalues and eigenvectors for the
fixed points  \\ 
of the replicator system with payoff matrix given by equation (4).$\text{%
\rule[-0.3in]{0in}{0.18in}}$%
\end{tabular}%
\end{center}

\end{document}